\documentclass[11pt,leqno,draft]{article}

\begin{document}

\title{Quasi-fractal sets in space}

\author{Stephen Semmes \\
        Rice University}

\date{}

\maketitle

        Let $a$ be a positive real number, $a < 1/2$.  A standard
construction of a self-similar Cantor set in the plane starts with the
unit square $[0, 1] \times [0, 1]$, replaces it with the four corner
squares with sidelength $a$, then replaces each of those squares with
their four corner squares of sidelength $a^2$, and so on.  At the
$n$th stage one has $4^n$ squares with sidelength $a^{-n}$, and the
resulting Cantor set has Hausdorff dimension $\log 4 / (-\log a)$.
The limiting case $a = 1/2$ simply reproduces the unit square, which
has Hausdorff dimension $2$.

        Suppose that we keep the boundaries of the squares at each
stage of the construction, to get a kind of quasi-fractal set
consisting of the Cantor set and a countable collection of line
segments.  The sum of the lengths of these line segments is finite
exactly when $a < 1/4$.  The Cantor set may be described as the
singular part of this quasi-fractal set, which is compact and
connected.

        Of course, one can consider similar constructions in higher
dimensions.  For the sake of simplicity, let us focus on connected
fractal sets in ${\bf R}^3$ with topological diimension $1$, for which
the corresponding quasi-fractal set is obtained by including a
countable collection of $2$-dimensional pieces.  As in the case of
higher-dimensional Sierpinski gaskets or Menger sponges, these
additional $2$-dimensional pieces could be triangles or squares.  Just
as the endpoints of the line segments were elements of the Cantor set
in the previous situation, the boundaries of these two-dimensional
pieces would be loops in the fractal.

        For that matter, one could do the same for connected fractals
in the plane like Sierpinski gaskets and carpets, where the additional
$2$-dimensional pieces are simply the bounded components of the
complement.  The resulting quasi-fractal would then be an ordinary
$2$-dimensional set, such as a triangle or a square.  In each of these
situations, the fractal set is the topological boundary of the rest,
which is disconnected but smooth.  One could also look at the bounded
components of the complement of the quasi-fractals in ${\bf R}^3$, but
we shall not pursue this here.

        This type of connected fractal set in the plane or space has a
lot of $1$-dimensional topological activity, which can be described in
terms of homotopy classes of continuous mappings to the circle.  In
particular, there are indices for such mappings associated to loops in
the fractal.  This suggests looking at Toeplitz operators and their
indices when they are Fredholm, using Bergman or Hardy spaces of
holomorphic functions on the two-dimensional smooth part.

        One already has complex analysis sitting on the plane, and for
the $2$-dimensional smooth pieces in space there are conformal
structures induced by the ambient Euclidean geometry.  Once
orientations are chosen, one gets complex structures, and one can talk
about Bergman and Hardy spaces.

        Different choices of orientations lead to different types of
Toeplitz operators.  It is natural to have a lot of different types of
Toeplitz operators in this situation, because of the complexity of the
fractals.


\begin{thebibliography}{21}


\bibitem {a} W.~Arveson, {\it A Short Course on Spectral Theory},
Springer-Verlag, 2002.

\bibitem {b-d-f-1} L.~Brown, R.~Douglas, and P.~Fillmore, {\it
Extensions of $C^*$-algebras, operators with compact self-commutators,
and $K$-homology}, Bulletin of the American Mathematical Society {\bf
79} (1973), 973--978.

\bibitem {b-d-f-2} L.~Brown, R.~Douglas, and P.~Fillmore, {\it Unitary
equivalence modulo the compact operators and extensions of
$C^*$-algebras}, in {\it Proceedings of a Conference on Operator
Theory}, 58--128, Lecture Notes in Mathematics {\bf 345},
Springer-Verlag, 1973.

\bibitem {b-d-f-3} L.~Brown, R.~Douglas, and P.~Fillmore, {\it
Extensions of $C^*$-algebras and $K$-homology}, Annals of Mathematics
(2) {\bf 105} (1977), 265--324.

\bibitem {c-w-1} R.~Coifman and G.~Weiss, {\it Analyse Harmonique
Non-Commutative sur Certains Espaces Homog\`enes}, Lecture Notes in
Mathematics {\bf 242}, Springer-Verlag, 1971.

\bibitem {c-w-2} R.~Coifman and G.~Weiss, {\it Extensions of Hardy
spaces and their use in analysis}, Bulletin of the American
Mathematical Society {\bf 83} (1977), 569--645.

\bibitem {cns} A.~Connes, {\it Noncommutative Geometry}, Academic
Press, 1994.

\bibitem {d1} R.~Douglas, {\it $C^*$-Algebra Extensions and
$K$-Homology}, Princeton University Press, 1980.

\bibitem {d2} R.~Douglas, {\it Banach Algebra Techniques in Operator
Theory}, 2nd edition, Springer-Verlag, 1998.

\bibitem {dr} P.~Duren, {\it Theory of $H^p$ Spaces}, Academic Press,
1970.

\bibitem {d-s} P.~Duren and A.~Schuster, {\it Bergman Spaces},
American Mathematical Society, 2004.

\bibitem {f} K.~Falconer, {\it The Geometry of Fractal Sets},
Cambridge University Press, 1986.

\bibitem {g} J.~Garnett, {\it Bounded Analytic Functions},
Springer-Verlag, 2007.

\bibitem {h-w} W.~Hurewicz and H.~Wallman, {\it Dimension Theory},
Princeton University Press, 1941.

\bibitem {k} J.~Kigami, {\it Analysis on Fractals}, Cambridge
University Press, 2001.

\bibitem {mat} P.~Mattila, {\it Geometry of Sets and Measures in
Euclidean Spaces}, Cambridge University Press, 1995.

\bibitem {sar} D.~Sarason, {\it Function Theory on the Unit Circle},
Virginia Polytechnic Institute and State University, 1978.

\bibitem {st1} E.~Stein, {\it Singular Integrals and Differentiability
Properties of Functions}, Princeton University Press, 1970.

\bibitem {st2} E.~Stein, {\it Harmonic Analysis: Real-variable
Methods, Orthogonality, and Oscillatory Integrals}, with the
assistance of T.~Murphy, Princeton University Press, 1993.

\bibitem {s-w} E.~Stein and G.~Weiss, {\it Introduction to Fourier
Analysis on Euclidean Spaces}, Princeton University Press, 1971.

\bibitem {str} R.~Strichartz, {\it Differential Equations on
Fractals}, Princeton University Press, 2006.



\end{thebibliography}
\end{document}